\documentclass{amsart}

\usepackage{pdfpages}
\usepackage[procnames]{listings}
\usepackage{color}
\usepackage{cite}
\usepackage{relsize}
\usepackage{amssymb}
\numberwithin{equation}{section}
\usepackage{amsmath}
\usepackage{mathtools}
\usepackage{amssymb}

\usepackage[applemac]{inputenc}
\usepackage[T1]{fontenc}

\title{ The least prime ideal in a given ideal class}
\author{Naser T. Sardari}
\date{\today}

	\newtheorem{thm}{Theorem}[section]
	
	\newtheorem{prop}[thm]{Proposition}

        \newtheorem{conj}[thm]{Conjecture}

	\newtheorem{rem}[thm]{Remark}

	\newtheorem{cor}[thm]{Corollary}
	
	\theoremstyle{defi}
	
	\newtheorem{cram}[thm]{Cram\'er's model}

	\theoremstyle{pf}

	\numberwithin{equation}{section}

\begin {document}
\maketitle
\begin{abstract}
 Let $K$ be a number field with the discriminant $D_K$ and the class number $h_{K}$, which has  bounded degree over $\mathbb{Q}$. By assuming GRH, we prove that every ideal class  of $K$ contains a prime ideal with norm less than $h_{K}^2\log(D_K)^{2}$ and also all but $o(h_K)$ of them have a  prime ideal with norm  less than $h_{K}\log(D_K)^{2+\epsilon}$.  For imaginary quadratic fields $K=\mathbb{Q}(\sqrt{D})$, by assuming  Conjecture~\ref{piarcor} (a weak version of the pair correlation conjecure), we improve our bounds by removing a factor of $\log(D)$ from our bounds and show that these bounds are optimal. 
%
\end{abstract}
\tableofcontents 
\section{Introduction}
\subsection{Motivation}
One application of our main result addresses  the optimal upper bound  on the least prime number represented by a binary quadratic form up to a power of the logarithm of its discriminant.  Giving sharp upper bound of this form on the least prime or the least integer representable by a sum of two squares is crucial in the analysis of the complexity of some algorithms in quantum compiling. In particular, Ross and Selinger's algorithm for the optimal navigation of $z$-axis rotations in $SU(2)$ by quantum gates \cite{RS} and its $p$-adic analogue for finding the shortest path between two diagonal vertices  of LPS Ramanujan graphs \cite{complexity}. In \cite{complexity}, we proved that these heuristic algorithms run in polynomial time under a Cram\'er type conjecture on the distribution of the inverse image of the integers representable as a sum of two squares by a binary quadratic from; see \cite[Conjecture 1.4.]{complexity}. In this paper, we fix the fundamental discriminant $D<0$, and prove that by assuming the generalized Riemann hypothesis a form of this Cram\'er type conjecture holds for all but $o(h(D))$ of binary quadratic forms with discriminant $D$ where $h(D)$ is the class number of $\mathbb{Q}(\sqrt{D})$. 
As a result, we prove that under GRH the above proposed algorithms give a probabilistic polynomial time algorithm for navigating $SU(2)$ and $PSL_2(\mathbb{Z}/q\mathbb{Z})$. We give a version of  our main theorem for the imaginary quadratic fields  $\mathbb{Q}(\sqrt{D})$ in what follows. Let $H_D$ be the  equivalence class of the binary quadratic forms with discriminant $D$ and    $\hat{H}_D$ denote the family of the unramified Hecke character associated to the  ideal class group of $\mathbb{Q}(\sqrt{D})$. Let   $R(D,X)$ denote the number of the classes of the binary quadratic forms of discriminant $D$ that does not represent any prime number $p$ where $ 1< p<X$. 
\begin{thm}\label{main} Let $X>D^{\alpha}$ for some $\alpha>0.$
Assume GRH holds for $L(s,\chi)$ where $\chi\in \hat{H}_D $ then 
$$R(D,X)  \ll \frac{h(D)^2 \log(|D|)^2}{X}.$$
Moreover, by assuming  Conjecture~\ref{piarcor}, we have 
$$R(D,X)  \ll \frac{h(D)^2 \log(|D|)}{X}.$$
\\
\end{thm}

%
%

The following conjecture is a weak version of the pair correlation conjecture for the family of the Hilbert characters associated to $\mathbb{Q}(\sqrt{D})$.  
\begin{conj} \label{piarcor} Let $w$ be a fixed smooth weight function where support $w\subset[1,2]$ and   $\hat{w}(s):=\int_{0}^{\infty}x^s w(x)\frac{dx}{x}$ be its Mellin transform. Assume GRH holds for $L(s,\chi)$ where $\chi\in \hat{H}(D)$ and $T>D^{\alpha}$ for some $\alpha>0$,  then 
\begin{equation}\label{corsum}
\frac{1}{h(D)}\sum_{\chi\in \hat{H}(D)}\sum_{\gamma_{\chi} ,\gamma_{\chi}^{\prime}}T^{i(\gamma_{\chi} -\gamma_{\chi}^{\prime})}\hat{w}(1/2+i\gamma_{\chi})\overline{\hat{w}(1/2+i\gamma_{\chi}^{\prime})}\ll \log(D),
\end{equation}
where $1/2+i\gamma_{\chi}$ and $1/2+i\gamma_{\chi}^{\prime}$ are zeros of  $L(s,\chi)$. Moreover,  $\ll$ is independent of  $\alpha$ and $D$ and only depends on  $w$.
\end{conj}
\begin{rem}
We expect the main term of \eqref{corsum} comes form the diagonal terms $\gamma_{\chi}=\gamma_{\chi}^{\prime}$. This is of order $\log(D)$ that is the number of the low lying zeros of $L(s,\chi)$ in a bounded interval.

\end{rem}

For $A\in H_D$, let $p_A$ be the smallest prime number representable by the class $A.$  
\begin{cor}\label{cor1}
Assume GRH holds for $L(s,\chi)$ where $\chi\in \hat{H}_D $ then for every $A\in H_D$ we have  
$$|p_A| \ll h(D)^2 \log(|D|)^2,$$
 and also for all $A\in H_D$ but $o(h(D))$ of them, we have 
$$|p_A| \ll h(D) \log(|D|)^{2+\epsilon},$$
Moreover by assuming  Conjecture~\ref{piarcor}, we have  $|p_A| \ll h(D)^2 \log(|D|)$ and $|p_A| \ll h(D) \log(|D|)^{1+\epsilon}$ for every and all but $o(h(D))$ of $A\in H_D.$ Moreover,  these bounds are optimal.
\end{cor}
\begin{rem}
 Littlewood \cite{Littlewood} proved that  under GRH $$L(1,\chi_{D})\ll \log\log(D)).$$ Hence by the class number formula, under GRH and Conjecture~\ref{piarcor}, almost all quadratic forms represent a prime number smaller than  $$\sqrt{|D|}\log(|D|)\log\log(|D|)f(D),$$ where
 $f(D)$ is any increasing function where $f(D)\to \infty$ as $D\to \infty.$
\\
\end{rem}

%
%
%

Recently, Thorner and Zaman esteblished  \cite[Theorem 1.2 ]{Thorner} the analogue of the Linnik's Thoerem for the binary quadratic forms of fixed discriminant. They proved that every positive definite binary quadratic form of discriminate $D$ represtes a prime number $p\ll D^{694}$ without any assumptions. Next,  we show that our bounds are optimal and compatible with the random model for the prime numbers known as  Cram\'er's model. 
\begin{prop}
Assume that  a positive proportion of the binary quadratic forms of discriminant $D<0$, represent a prime number less than $X$ then 
$$h(D)\log{D}\ll X.$$
\end{prop}
\begin{proof}
 Let $r(n,D)$ denote the sum of the representation of $n$ by all the classes of binary quadratic forms of discriminant $D$
$$r(n,D)=\sum_{Q\in H_D}r(n,Q).$$
By the classical formula due to Dirichlet we have
$$r(n,D)=w_D\sum_{d|n}\chi_{D}(d),$$  
where,
$$w_D=\begin{cases} 6 \text{ if }D=-3 \\  4 \text{ if }D=-4 \\   2 \text{ if } D<-4.  \end{cases}$$
This means that the multiplicity of representing a prime number $p$ by all the binary quadratic forms of a fixed negative discriminant $D<-4$ is bounded by $4$
\begin{equation}\label{primebound}
r(p,D) \leq 4.
\end{equation}
Assume that a positive proportion of binary quadratic forms represent a prime number smaller than $X$. Let $N(X,D)$ denote the number of pairs $(p,Q)$ such that $p$ is a number less than $X$ and $Q\in H_D$ is a representative of a binary quadratic form of discriminant $D$. We proceed by giving a double counting formula for $N(X,D)$. By our assumption a positive proportion of binary quadratic forms of discriminant $D$ represent a prime number less than $X$, then
\begin{equation}\label{low}h(D)\ll N(X,D).\end{equation}
On the other hand,
$$N(X,D)=\sum_{p<X} r(p,D).$$
By inequality \eqref{primebound},
$$N(X,D)\leq 4 \pi(X).$$
By the above inequality and inequality \eqref{low}, we obtain
$$h(D)\ll \pi(X).$$
By Siegel's lower bound $D^{1/2-\varepsilon}\ll h(D)$, the above inequality is equivalent to 
$$h(D)\log(D)\ll X.$$
This completes the proof of our proposition.
\end{proof}
 We cite the following formulation of the  Cram\'er model from~\cite{Soundist}. 
\begin{cram}
The primes behave like independent random variables $X(n)$ $ (n \geq 3)$ with $X(n) = 1$ (the number $n$ is `prime') with probability $1/\log n$, and $X(n) = 0$ (the number $n$ is `composite') with probability $1-1/\log n.$
\end{cram}

 Note that each class of the integral binary quadratic forms is associated to  a Heegner point in $SL_2(\mathbb{Z}) \backslash \mathbb{H}$. By the equdistribution of Heegner points in $SL_2(\mathbb{Z}) \backslash \mathbb{H}$, it follows that almost all classes of the integral quadratic forms has a representative $Q(x,y):=Ax^2+Bxy+Cy^2$ such that the coefficients of $Q(x,y)$ are bounded by any function growing faster than $\sqrt{D}$:
$$\max(|A|,|B|,|C|)<\sqrt{D}\psi(D),$$
for any function $\psi(D)$ defined on integers such that $\psi(D)\to \infty$ as $D \to \infty$. We show this claim in what follows. We consider the set of representative of the Heegner points inside the Gauss fundamental domain of $SL_2(\mathbb{Z}) \backslash \mathbb{H}$ and denote them by $z_{\alpha}$ for $\alpha \in H(D)$. They are associated to the roots of a representative of a binary quadratic form in the ideal class group. By the equidistribution of Heegner points in $SL_2(\mathbb{Z}) \backslash \mathbb{H}$ and the fact that the volume of the  Gauss fundamental domain decay with rate $y^{-1}$ near the cusp, it follows that for almost all $\alpha \in H(D)$ if $z_{\alpha}=a+ib$ is the Heegner point inside the Gauss fundamental domain  associated to $\alpha$  then
\begin{equation}
\begin{split}
\label{bineq}
|a| \leq 1/2,
\\
\sqrt{3}/2 \leq  b\leq \psi(D),
\end{split}
\end{equation}
where $\psi(D)$ is any function such that $\psi(D)\to \infty$ as $D\to \infty.$ 
Let $Q_{\alpha}(x,y):=Ax^2+Bxy+Cy^2$ be the quadratic forms associated to $\alpha \in H(D)$ that has $z_{\alpha}$ as its root. Then 
$$z_{\alpha}=\frac{-B\pm i\sqrt{D}}{2A},$$
where $a=\frac{-B}{2A}$ and $b=\frac{\sqrt{D}}{2A}$.
By inequality~\eqref{bineq}, we have
\begin{equation}
\begin{split}
|B|\leq |A|,
\\
\frac{\sqrt{D}}{2\psi(D)} \leq A < \sqrt{D}.
\end{split}
\end{equation}
By the above inequalities and $D=B^2-4AC$, it follows that
\begin{equation}\label{hbd}\max(|A|,|B|,|C|)<\sqrt{D}\psi(D).\end{equation}
This concludes our claim. Next, we give a heuristic upper bound on the size of the smallest prime number represented by a binary quadratic forms of discriminant $D$ that satisfies \eqref{hbd}.  Since $D$ is square-free, there is no local restriction for representing prime numbers.   So, by the Cram\'er's model and consideration of the Hardy-Littlewood local measures,  we expect that for a positive proportion of the classes of the binary quadratic forms $Q$ there exists an integral point $(a,b) \in \mathbb{Z}^2$ such that $|(a,b)|^2<L(1,\chi_{D})\log(D)$ and  $Q(a,b)$ is a prime number. We have 
\begin{equation}\label{Cramer}
\begin{split}
Q(a,b)&=Aa^2+Bab+Cb^2
\\
&\leq \max(|A|,|B|,|C|)|(a,b)|^2
\\
&\ll\sqrt{D}L(1,\chi_{D})\psi(D)\log(D).
\end{split}
\end{equation}
We may take $\psi(D)$ to be any constant in the above estimate. Therefore, we expect that almost all quadratic forms of discriminant $D$ represent a prime number less than $h(D)\log(D)^{1+\epsilon}.$  In \cite{least}, we proved unconditionally  that a positive proportion of the  binary quadratic forms of discriminant  $D<0$  represent a prime number smaller than any  fixed scalar multiple of the optimal size  $h(D)\log(D)$, where $h(D)$ is the class number of $\mathbb{Q}(\sqrt{-D})$. 
By a similar analysis, we expect that almost all ideal class of the a bounded degree number field $K$ with discriminant  $D_K$ contain  a prime number less than 
$h_K\log(D_K)^{1+\epsilon}.$

\subsection{Repulsion of the prime ideals near the cusp }As we noted above, based on the Cram\'er's model we expect that the split prime numbers  randomly distributed among the ideal classes of $\mathbb{Q}(\sqrt{D})$. By comparing with the coupon collector's problem,   we may expect  that every ideal class contain a prime ideal of size $h(D)D^{\epsilon}$. Note that Cram\'er conjecture states that every short interval of size $\log(X)^{2+\epsilon}$ contains a prime number. By Linnik's conjecture, every congruence class modulo $q$ contains a prime number less than $q^{1+\epsilon}.$ This shows that small prime numbers covers all the short interval and congruence classes. However,  small primes of size $h(D)D^{\epsilon}$  are not covering all  the classes of binary quadratic forms. For example, the least prime number represented by $Q(x,y)=Dx^2+y^2$  is bigger than $D$ compared to $\sqrt{D}\log(D)^{2+\varepsilon}$ that is the upper bound for almost all binary quadratic forms  under GRH. This feature is different from the analogues conjectures for the size of the least prime number in a given congruence classes modulo an integer (Linnik's conjecture) and the distribution of prime numbers in short intervals (Cram\'er's conjecture). We call this new feature the repulsion of small primes by the cusp. In fact, the binary quadratic forms with the associated Heegner point near the cusp repels prime numbers. This can be seen in equation \eqref{Cramer}, where $\max(|A|,|B|,|C|)|$ could be as large as $D$ near the cusp whereas  for a typical binary quadratic form it is bounded by $D^{1/2+\epsilon}.$ This shows that  both bounds for the almost all and every class of binary quadratic forms in Corollary~\ref{cor1} are sharp.

Finally, by the one to one correspondence between the classes  of the binary quadratic forms with discriminant $D$ and the ideal classes of $\mathbb{Q}(\sqrt{D})$, we prove a Minkowski's bound on the least prime ideal  in a given ideal class of $\mathbb{Q}(\sqrt{D})$. More precisely, let $D<0$ be square-free and  $D\equiv 1 \mod 4$. Let $H(D)$ denote the ideal class group of $\mathbb{Q}(\sqrt{D})$ and $N_{\mathbb{Q}(\sqrt{D})} (x+y\sqrt{D})=x^2-Dy^2$ be the norm of the imaginary quadratic field $\mathbb{Q}(\sqrt{D}).$ Given an integral ideal $I\subset \mathcal{O}_{\mathbb{Q}(\sqrt{D})}$, let $q_{I}(x,y)$ be the following class of the integral binary quadratic form defined up to the action of $SL_2(\mathbb{Z})$ 
\begin{equation}\label{classis}
q_{I}(x,y):=\frac{N_{\mathbb{Q}(\sqrt{D})}(x\alpha+y\beta)}{N_{\mathbb{Q}(\sqrt{D})}(I)}\in \mathbb{Z},
\end{equation}
where $x,y\in \mathbb{Z}$, and $ I \cong \langle \alpha,\beta \rangle_{\mathbb{Z}}$ identifies the integral ideal $ I $ with $ \mathbb{Z}^2$.   It follows that $q_{I}$ only depends on the ideal class $[I]\in H(D).$ This gives an isomorphism between $H(D)$ and $H_D$ (the class of binary quadratic forms of discriminant $D$). 
%
%
 In fact in our main theorem, we  prove this type of generalized Minkowski's bounds for every number fields. More precisely, let $K$ be a number field with the discriminant $D_K$ and the class number $h_{K}$, which has  bounded degree over $\mathbb{Q}$.  Let $H_{K}$ be the  ideal class group of $K$ and $\hat{H}_K$ be the character group of the ideal class group $H_K$ that is also identified with the the maximal unramified characters (Hilbert characters) of the number field $K$;  see \cite{hecke}. 
   Minkowski proved that every $A\in H_K$ contains an integral ideal of norm 
$
O(|D_K|^{1/2})
$ where the implicit constant in $O$ only depends on the degree of $K$ over $\mathbb{Q}.$
 Let $R(K,T)$ denote the number of $A\in H_K$ that does not contain  a prime ideal $\mathfrak{p}$ where $ 1< N_K(\mathfrak{p})<T$. 
We prove the following generalization of the Minkowski's bound for any  number field $K$  with bounded degree over $\mathbb{Q}$.   
\begin{thm}\label{genmisk}
Let $K$ be as above.  Assume GRH holds for $L(s,\chi)$ where $\chi\in \hat{H}_K $ then 
$$R(K,T)  \ll \frac{h_K^2 \log(|D_K|)^2}{T}.$$
\\
\end{thm}
For $A\in H_K$, let $\mathfrak{p}_A$ be the smallest prime ideal inside  the class $A.$ The following is an immediate corollary of Theorem~\ref{genmisk}.
\begin{cor}\label{leastmin}
Assume GRH holds for $L(s,\chi)$ where $\chi\in \hat{H}_K $ then for every $A\in H_K$ we have  
$$|\mathfrak{p}_A| \ll h_K^2 \log(|D_K|)^2,$$
 and also for all $A\in H_K$ but $o(h_K)$ of them, we have 
$$|\mathfrak{p}_A| \ll h_K \log(|D_K|)^{2+\epsilon}.$$
\end{cor}  
\begin{proof}
 By choosing $T\gg h_K^2 \log(D_K)^{2}$ in  Theorem~\ref{genmisk}, we deduce that  
 $
 R(K,T)=0. 
 $ This shows that for every $A\in H_K$ we have  
$|\mathfrak{p}_A| \ll h_K^2 \log(|D_K|)^2.$ For the next inequality, let $X=h_K \log(|D_K|)^{2+\epsilon}$ and apply Theorem~\ref{genmisk} to obtain 
$$
R(K,X)  \ll \frac{h_K^2 \log(|D_K|)^2}{X}\leq \frac{h_K}{ \log(|D_K|)^{\epsilon}}=o(h_K),
$$
which concludes Corollary~\ref{leastmin}.
\end{proof}

\subsection{Method and Statement of results}

 Our strategy of the proof is based on Selberg's method \cite{Selberg} for investigating the distribution of primes in short intervals. By Cram\'er's conjecture,  every short interval of size $\log(X)^{2+\varepsilon}$  inside the large  interval $[X,2X]$ contains a prime number. Selberg proved that by assuming the Riemann hypothesis  almost all short intervals  of size $\log(X)^{2+\varepsilon}$ in the large  interval $[X,2X]$ contain a prime number. By assuming the Riemann hypothesis, Selberg gave a sharp upper bound up to a power of $\log(X)$ on the variance of the number of prime numbers inside a short interval.  The result follows from  Chebyshev's inequality and this sharp upper bound on the variance.  Heath-Brown showed that this extra $\log(X)$ can be removed by assuming  Montgomery's pair correlation conjecture \cite[Corollary 4]{Heatpair}. Namely, by assuming Riemann hypothesis and the pair correlation conjecture for the Riemann zeta function almost all short intervals  of size $\log(X)^{1+\varepsilon}$ in the large  interval $[X,2X]$ contain a prime number. Similarly, by assuming GRH for the Dirichlet L-function $L(s,\chi)$ with conductor $q$, Montgomery \cite[Theorem 17.1]{Mon}, gave a sharp upper bound up to a factor of $\log(q)$ on the variance of the number of prime numbers distributed in difference congruence classes modulo $q.$ It follows that the least prime in almost all congruence classes modulo  $q$ is less than $\varphi(q)\log(q)^{2+\varepsilon}$. Our result follows from the analogues upper bound on the variance of the number of prime ideals in ideals classes.  By assuming, Conjecture~\ref{piarcor} on the distribution of the low lying zeros of the family of Hecke L-function associated to $\mathbb{Q}(\sqrt{D})$ we  also save a factor of $\log(D)$. 
\subsection*{Acknowledgements}\noindent
I would like to thank Professor Roger Heath-Brown my mentor at MSRI for several insightful and inspiring conversations during the Spring~2017 Analytic Number Theory program at MSRI. This material is based upon work supported by the National Science Foundation under Grant No. DMS-1440140 while the author  was in residence at the Mathematical Sciences Research Institute in Berkeley, California, during the Spring 2017 semester.

\section{Hecke L-functions}
\subsection{Functional equation}
In this section, we review the basic properties of the Hecke L-functions with Gr\"ossencharaktere, known as Hecke characters.  In particular, the functional equation of the unramified characters (Hilbert characters) associated to the ideal class group of a number field $K$ over $\mathbb{Q}$ and Weil's explicit formula. We begin by introducing some notations. Recall that $K$ is a number field of the bounded degree $n=r_1+2r_2$ over $\mathbb{Q}$  with the  discriminant $D_K$ and the class number $h_{K}$. We write $N_{K}$ for the norm forms defined on $K.$  Let $H_{K}$ be the  ideal class group of $K$ and $\hat{H}_K$ be the character group of the ideal class group $H_K$ that is also identified with the the maximal unramified characters (Hilbert characters) of the number field $K$;  see \cite{hecke}. 
 Assume that $\chi\in \hat{H}_K$ is a Hilbert character. Let 
\begin{equation}
L(s,\chi):=\sum_{\mathfrak{a}} \chi(\mathfrak{a})N_{K}(\mathfrak{a})^{-s},
\end{equation}
for $\Re(s)>1$ where the sum is over integral ideals $\mathfrak{a}$ of the ring of integers $\mathcal{O}_K.$ By the unique factorization of the ideals into the prime ideals, $L(s,\chi)$ has an Euler product for $\Re(s)>1$
\begin{equation}
L(s,\chi)=\prod_{\mathfrak{p}} \big(1-\chi(\mathfrak{p})N_{K}(\mathfrak{p})^{-s}   \big)^{-1}.
\end{equation}
Hecke proved that $L(s,\chi)$ satisfies the following functional equation; see \cite[equation (1.2.2)]{MR1031335}
\begin{equation}
L(s,\chi)=w(\chi)A^{1-2s}G(1-s,\chi_{\infty})L(1-s,\bar{\chi}),
\end{equation}
where $|w(\chi)|=1$, $A^2=|D_k|\pi^{-n}2^{-r_2}$ and $G(1-s,\chi_{\infty})$ is the contribution of the gamma factor at $\infty.$ Since we assumed that $\chi$ is a Hilbert character then $\chi$ has trivial  weight and $\chi_{\infty}=1$ and we have 
\begin{equation}\label{gmf}
G(s,\chi_{\infty}):= \Big(\frac{\Gamma(\frac{1}{2}s )}{\Gamma(\frac{1}{2}(1-s) )}\Big)^{r_1} \Big(\frac{\Gamma(s)}{\Gamma(1-s)}   \Big)^{r_2}.
\end{equation}
From the above functional equation, we obtain the following asymptotic formula for the number of the zeros of $L(s,\chi)$; see \cite[Proposition 2.3.1] {Kowlaski}.  
\begin{prop}\label{zeronum} We have
\begin{equation}
N(\chi,T)=\frac{T}{2\pi} \log\big(\frac{T^n|D_K|}{(2\pi e)^n}   \big)+O\big(\log (T^n|D_K|)   \big),
\end{equation}
for $T\geq 2$, the estimate being uniform in all parameters.

\end{prop}

\subsection{Weil's explicit formula}

We cite the following version of the explicit formula for the  Hecke  L-function with Gr\"ossencharaktere. \cite[Proposition 2.3.5]{Kowlaski}
\begin{prop}\label{weil} 
Let $w: (0,\infty)\to \mathbb{C}$ be a $C^{\infty}$ function with compact support, let 
$$
\hat{w}(s)=\int_{0}^{+\infty}w(x)x^{s-1}dx,
$$
be its Mellin transform, which is entire and decays rapidly in vertical strips. Let 
$$
\varphi(x)=\frac{1}{x}w(\frac{1}{x}).
$$
Assume that $\chi\in \hat{H}_K$ is a Hilbert character. Then we have
\begin{equation}\label{explicit}
\begin{split}
\sum_{\mathfrak{a}} \Lambda(\mathfrak{a})\big(\chi(\mathfrak{a})w(N_K(\mathfrak{a}))+\bar{\chi}(\mathfrak{a})\varphi(N_K(\mathfrak{a}))   \big)=(\log|D_K|)w(1)+\delta(\chi)\int_0^{\infty}w(x)dx
\\
-\sum_{\rho} \hat{w}(\rho)+\frac{1}{2\pi i } \int_{\Re(s)=-1/2} \frac{G^{\prime}(s,\chi_{\infty})}{G(s,\chi_{\infty})} \hat{w_T}(s)ds +\text{Res}_{s=0}\big(-\frac{L^{\prime}}{L}(\chi,s)\hat{w}(s)   \big), 
\end{split}
\end{equation}
where $\delta(\chi)=1$ if $\chi$ is trivial and is $=0$ otherwise, and the sum over zeros includes multiplicity. 
\end{prop}

\section{Proof of Theorem~\ref{genmisk}}
\begin{proof}
  Let $w(x)$ be a fixed positive real valued smooth compactly supported function supported  on the interval $[1,2]$ such that $\int_{1}^{2}w(x)dx=1$ and let $w_T(x):=w(x/T)$. For an ideal class  $A\in H_K,$ we define
 $$\psi_{A}(w_T):= \sum_{\mathfrak{n}\in A}\Lambda(\mathfrak{n})w_T(N_{K}(\mathfrak{n})), $$
 where  $\Lambda(\mathfrak{n})$ is the generalized von Mangoldt function defined by
 \begin{equation}
 \Lambda(\mathfrak{n})=\begin{cases}
 \log N_{K}(\mathfrak{p}), \quad   &\mathfrak{a}=\mathfrak{p}^ k \text{ for some prime ideal } \mathfrak{p}\subset \mathcal{O}_K,
 \\
 0,  & \text{otherwise.}
 \end{cases}
 \end{equation}
Moreover, for every Hilbert character $\chi \in \hat{H}_K$, we define
\begin{equation}\label{Fourier}\psi_{\chi}(w_T):= \sum_{\mathfrak{n}} \chi(\mathfrak{n})\Lambda(\mathfrak{n})w_T(N_{K}(\mathfrak{n})). \end{equation}
Finally, we define  
$$\psi(w_T):=\sum_{\mathfrak{n}} \Lambda(\mathfrak{n})w_T(N_{K}(\mathfrak{n})), $$
which is associated to the trivial character on $H_K$ and also  the sum of  $\psi_{A}(w_T)$ over $A\in H_K$.   It follows from the explicit formula in Proposition~\ref{weil} that by assuming GRH for the Dedekind zeta function $\zeta_K$; see~\cite[Theorem 1.1]{LO}
\begin{equation}
\big| \psi(w_T)-T  \big| \leq C_0T^{1/2}\log(D_KT^{n}).
\end{equation}
In particular, if $T\gg \log(D_K)^{2+\epsilon}$ then 
$$
\psi(w_T)=T(1+o(1)).
$$
 By the orthogonality of characters, we express $\psi_{A}(w_T)$ in terms of $\psi_{\chi}(w_T)$
\begin{equation}\label{Fourier1}\psi_{A}(w_T):= \frac{1}{h_K}\sum_{\chi\in \hat{H}_K}\bar{\chi}(A)\psi_{\chi}(w_T). \end{equation}
We define $\text{Var}(K,w_T)$ as
\begin{equation}\label{Varr} \text{Var}(K,w_T):=\sum_{A\in H_K} |\psi_{A}(w_T)-\frac{1}{h_K}\psi(w_T)|^2. \end{equation}
%
We use equation~(\ref{Fourier1}) and obtain
\begin{equation*}
\begin{split}
 \text{Var}(K,w_T)&=\sum_{A\in H_K} \Big|\frac{1}{h_K} \sum_{\chi\neq id  }\bar{\chi}(A) \psi_{\chi}(w_T) \Big|^2 
 \\
 &=\frac{1}{h_K^2} \sum_{\chi_1,\chi_2\neq id }   \psi_{\chi_1}(w_T) \bar{\psi}_{\chi_2}(w_T) \sum_{A\in H_K} \bar{\chi}_1\chi_2(A).
  \end{split}
 \end{equation*}
 By the orthogonality of characters, 
$$  \sum_{A\in H_K} \bar{\chi}_1\chi_2(A)=\delta(\chi_1,\chi_2):=\begin{cases}1 &\text{ if } \chi_1=\chi_2, \\ 0  &\text{ Otherwise.}  \end{cases}$$
 Hence, only the diagonal terms $\chi_1=\chi_2$ contribute to $ \text{Var}(K,w_T)$ and we obtain
\begin{equation}\label{Var}  \text{Var}(K,w_T)= \frac{1}{h_K} \sum_{\chi\neq id} |\psi_{\chi}(w_T) |^2. \end{equation}
By applying the explicit formula \eqref{explicit}, we relate $|\psi_{\chi}(w_T) |$  to the zeros of $L(s,\chi)$
\begin{equation}\label{mellin}
\begin{split}
\psi_{\chi}(w_T)=-\psi_{\bar{\chi}}(\varphi_T) +(\log|D_K|)w_T(1)+\delta(\chi)\int_0^{\infty}w_T(x)dx -\sum_{\rho} \hat{w_T}(\rho)+\\
+\frac{1}{2\pi i } \int_{\Re(s)=-1/2} \frac{G^{\prime}(s,\chi_{\infty})}{G(s,\chi_{\infty})} \hat{w_T}(s)ds +\text{Res}_{s=0}\big(-\frac{L^{\prime}}{L}(\chi,s)\hat{w_T}(s)   \big), 
\end{split}
\end{equation}
where $ \hat{w_T}$ is the Mellin transform of $w_T$ defined as 
$
\hat{w_T}(s):=\int_{0}^{\infty}x^s w_T(x)\frac{dx}{x},
$
and 
$
\varphi_T(x):=\frac{1}{x}w_T(\frac{1}{x}).
$
Next, we analyze the terms on the right hand side of \eqref{mellin}. Since our weight function $w$ is supported on $[1,2]$, then $w_T$ is supported on $[T,2T]$ and also $\varphi_T$ is supported on $[1/2T, 1/T]$. Hence the first and the second terms vanishes in  \eqref{mellin} 

\begin{equation*}
\begin{split}
\psi_{\bar{\chi}}(\varphi_T)&= \sum_{\mathfrak{n}} \chi(\mathfrak{n})\Lambda(\mathfrak{n})\varphi_T(N_{K}(\mathfrak{n}))= \sum_{\mathfrak{n}} \chi(\mathfrak{n})\Lambda(\mathfrak{n})\frac{1}{N_{K}(\mathfrak{n})} w(\frac{1}{TN_{K}(\mathfrak{n})})=0.
\end{split}
\end{equation*}
Since $\int w(x)dx=1$, then 
\begin{equation}\label{delta}
\delta(\chi)\int_0^{\infty}w_T(x)dx= \delta(\chi)T.
\end{equation}
%
Recall that $w_T(x):=w(x/T)$ then
$
\hat{w_T}(s)=T^{s}\hat{w}(s).
$
Moreover,  by GRH all the non-trivial  zeros $\rho$ of $L(s,\chi)$ are on the critical line $Re(z)=1/2$. Hence, 
 \begin{equation}
 \begin{split}
\sum_{\rho} \hat{w_T}(\rho)= \sum_{\rho=1/2+i\gamma}T^{1/2+i\gamma} \hat{w}(1/2+i\gamma)
 \\
\leq T^{1/2}\Big(\sum_{\gamma} |\hat{w}(1/2+i\gamma)|    \Big),
\end{split} 
\end{equation}
 where $\rho=1/2+i\gamma$ is a zero of $L(s,\chi)$. Note that the smooth sum $\sum_{\gamma} |\tilde{w}(1/2+i\gamma)| $ over the zeros of $L(s,\chi)$ is localized on a fixed interval on the critical line   $Re(z)=1/2.$  Therefore it is bounded  (up to a constant depending on $w$) by  the number of zeros of the $L(s,\chi)$ in a fixed interval. By Proposition~\ref{zeronum}, the number of low lying zeros in a fixed interval is bounded by the logarithm of the conductor of $L(s,\chi)$. Therefore,
 \begin{equation}\label{sum}
 \Big|\sum_{\rho} \hat{w_T}(\rho)   \Big| \ll T^{1/2}\log (|D_K|),
 \end{equation}
 where the implicit constant in $\ll$ only depends on the smooth weight function $w$ and $n=[K:\mathbb{Q}].$ Finally, we analyze the contribution of the gamma factors in the Weil's explicit formula.  We have
 \begin{equation*}
 \begin{split}
  \int_{\Re(s)=-1/2} \frac{G^{\prime}(s,\chi_{\infty})}{G(s,\chi_{\infty})} \hat{w_T}(s)ds&\leq T^{-1/2} \int_{\Re(s)=-1/2} \big| \frac{G^{\prime}(s,\chi_{\infty})}{G(s,\chi_{\infty})} \hat{w}(s) \big|ds.
  \end{split}
 \end{equation*}
 By equation~\eqref{gmf}, we have 
 $$
\frac{G^{\prime}(s,\chi_{\infty})}{G(s,\chi_{\infty})}=r_1\big(\frac{\Gamma^{\prime}(s/2)}{\Gamma(s/2)} -\frac{\Gamma^{\prime}((1-s)/2)}{\Gamma((1-s)/2)}  \big)+r_2\big(\frac{\Gamma^{\prime}(s)}{\Gamma(s)} -\frac{\Gamma^{\prime}(1-s)}{\Gamma(1-s)}  \big)
 $$
By Stirling's formula, it follows that
  $$
|\frac{G^{\prime}(-1/2+it,\chi_{\infty})}{G(-1/2+it,\chi_{\infty})}| \ll n\log(1+|t|),
  $$
  where the implied constant is absolute. Hence, 
  \begin{equation}\label{gammafact}
   \big|\int_{\Re(s)=-1/2} \frac{G^{\prime}(s,\chi_{\infty})}{G(s,\chi_{\infty})} \hat{w_T}(s)ds \big|  = O(T^{-1/2}),
  \end{equation}
  where the implied constant in $O$ depends only on $n$ and the smooth function $w.$ Finally, we have 
  \begin{equation}\label{resi}
  \begin{split}
 \text{Res}_{s=0}\big(-\frac{L^{\prime}}{L}(\chi,s)\hat{w_T}(s)   \big)&=   \text{Res}_{s=0}\big(-\frac{L^{\prime}}{L}(\chi,s)T^s\hat{w}(s)   \big)
 \\
 &=\begin{cases}
 O(\log(T)^{r_1+r_2-1})  &\text{ if } \chi \text{ is nontrivial  }
 \\
 O(\log(T)^{r_1+r_2-2} )  &\text{ if } \chi \text{ is the trivial character },
 \end{cases}
 \end{split}
  \end{equation}
  where the implicit constant in $O$ only depends on $w.$ By combining upper bounds \eqref{gammafact}, \eqref{sum}, \eqref{gammafact} and \eqref{resi}, we obtain 
  \begin{equation}\label{finall}
|\psi_{\chi}(w_T)-\delta{\chi}T|= O(T^{1/2}\log (|D_K|)).
  \end{equation}
 By using the above upper bound in equation~(\ref{Var}), we obtain
  \begin{equation}\label{upper}
  \text{Var}(K,w_T) \ll T \log(D_K)^2.
  \end{equation}
 Recall from equation~(\ref{Varr}),
 \begin{equation} 
 \begin{split}
  \text{Var}(K,w_T)&=\sum_{Q\in H_D} |\psi_{Q}(w_T)-\frac{\psi(w_T)}{h_K}|^2 
  \\
  &\gg R(K,T) |\frac{\psi(w_T)}{h_K}|^2
  \end{split}
  \end{equation}
By equation~\eqref{finall}, $|\psi(w_T)|\sim T$. Therefore,
  \begin{equation}\label{lower}  
   \text{Var}(K,w_T)\gg R(K,T)\frac{T^2}{h_K^2}.
  \end{equation}
 By inequalities~(\ref{upper}) and (\ref{lower}), we obtain
 \begin{equation}
 R(K,T) \frac{T^2}{h_K^2} \ll T \log(D_K)^2.
 \end{equation}
Hence,
 $$R(K,T)  \ll \frac{h_K^2 \log(D_K)^2}{T}.$$
 This concludes the proof of our theorem. 

 \end{proof}
\section{Proof of Theorem~\ref{main} } 
\begin{proof} The first inequality is a consequence of Theorem~\ref{genmisk}. For the second inequality,  we follow a similar method as in the proof of Theorem~\ref{genmisk} and define the variance $ \text{Var}(\mathbb{Q}(\sqrt{D}),w_T)$ as in \eqref{Varr}.  By identity~\eqref{Var},  we obtain
\begin{equation}\label{varD} \text{Var}(\mathbb{Q}(\sqrt{D}),w_T)= \frac{1}{h(D)} \sum_{\chi\neq id} |\psi_{\chi}(w_T) |^2.\end{equation}
We apply Weil's explicit formula as in \eqref{mellin} and obtain 
\begin{equation*}
\begin{split}
\psi_{\chi}(w_T)=-\psi_{\bar{\chi}}(\varphi_T) +(\log|D|)w_T(1)+\delta(\chi)\int_0^{\infty}w_T(x)dx -\sum_{\rho} \hat{w_T}(\rho)+\\
+\frac{1}{2\pi i } \int_{\Re(s)=-1/2} \frac{G^{\prime}(s,\chi_{\infty})}{G(s,\chi_{\infty})} \hat{w_T}(s)ds +\text{Res}_{s=0}\big(-\frac{L^{\prime}}{L}(\chi,s)\hat{w_T}(s)   \big).
\end{split}
\end{equation*}
As before, the first and the second terms vanishes in the above formula and by applying the bounds \eqref{delta}, \eqref{gammafact} and \eqref{resi}, we have  
$$
\psi_{\chi}(w_T)- \delta(\chi)T=\sum_{\gamma_{\chi}}T^{1/2+i\gamma_{\chi}} \hat{w}(1/2+i\gamma_{\chi})+  O(\log(T)^{r_1+r_2-1}). 
$$
We square the above identity and obtain
\begin{equation*}\label{finaleq}
|\psi_{\chi}(w_T)- \delta(\chi)T|^2= \sum_{\gamma_{\chi} ,\gamma_{\chi}^{\prime}}T^{1+i(\gamma_{\chi} -\gamma_{\chi}^{\prime})}\hat{w}(1/2+i\gamma_{\chi})\hat{w}(1/2+i\gamma_{\chi}^{\prime})+O(T^{1/2}\log(T)^{r_1+r_2}).
\end{equation*}
By averaging  the above identity over $\chi\neq id$, we obtain
\begin{equation*}
\text{Var}(\mathbb{Q}(\sqrt{D}),w_T)=  \frac{1}{h(D)} \sum_{\chi\neq id}  \sum_{\gamma_{\chi} ,\gamma_{\chi}^{\prime}}T^{1+i(\gamma_{\chi} -\gamma_{\chi}^{\prime})}\hat{w}(1/2+i\gamma_{\chi})\hat{w}(1/2+i\gamma_{\chi}^{\prime})+O(T^{1/2}\log(T)^{r_1+r_2}).
\end{equation*}
By Conjecture~\ref{piarcor}, we obtain 
 $$
 \text{Var}(\mathbb{Q}(\sqrt{D}),w_T)\ll T\log(D).
 $$
 The above inequality save a factor of $\log(D)$ comparing to the inequality \eqref{upper}. This concludes  the second inequality in  Theorem~\ref{main}.%
%
\end{proof}
\bibliographystyle{alpha}
\bibliography{Minkowski}

\begin{thebibliography}{{Sar}18}

\bibitem[Duk89]{MR1031335}
W.~Duke.
\newblock Some problems in multidimensional analytic number theory.
\newblock {\em Acta Arith.}, 52(3):203--228, 1989.

\bibitem[HB82]{Heatpair}
D.~R. Heath-Brown.
\newblock Gaps between primes, and the pair correlation of zeros of the zeta
  function.
\newblock {\em Acta Arith.}, 41(1):85--99, 1982.

\bibitem[Hec20]{hecke}
E.~Hecke.
\newblock Eine neue {A}rt von {Z}etafunktionen und ihre {B}eziehungen zur
  {V}erteilung der {P}rimzahlen.
\newblock {\em Math. Z.}, 6(1-2):11--51, 1920.

\bibitem[ILS00]{density1}
Henryk Iwaniec, Wenzhi Luo, and Peter Sarnak.
\newblock Low lying zeros of families of {$L$}-functions.
\newblock {\em Inst. Hautes \'Etudes Sci. Publ. Math.}, (91):55--131 (2001),
  2000.

\bibitem[Kow03]{Kowlaski}
E.~Kowalski.
\newblock Elementary theory of {$L$}-functions. {II}.
\newblock In {\em An introduction to the {L}anglands program ({J}erusalem,
  2001)}, pages 21--37. Birkh\"auser Boston, Boston, MA, 2003.

\bibitem[KS99]{density}
Nicholas~M. Katz and Peter Sarnak.
\newblock {\em Random matrices, {F}robenius eigenvalues, and monodromy},
  volume~45 of {\em American Mathematical Society Colloquium Publications}.
\newblock American Mathematical Society, Providence, RI, 1999.

\bibitem[Lit28]{Littlewood}
J.~E. Littlewood.
\newblock On the class-number of the corpus p(??k).
\newblock {\em Proceedings of the London Mathematical Society},
  s2-27(1):358--372, 1928.

\bibitem[LO77]{LO}
J.~C. Lagarias and A.~M. Odlyzko.
\newblock Effective versions of the {C}hebotarev density theorem.
\newblock pages 409--464, 1977.

\bibitem[Mon71]{Mon}
Hugh~L. Montgomery.
\newblock {\em Topics in multiplicative number theory}.
\newblock Lecture Notes in Mathematics, Vol. 227. Springer-Verlag, Berlin-New
  York, 1971.

\bibitem[RS14]{RS}
N.~J. {Ross} and P.~{Selinger}.
\newblock {Optimal ancilla-free Clifford+T approximation of z-rotations}.
\newblock {\em ArXiv e-prints}, March 2014.

\bibitem[{Sar}17]{complexity}
N.~T {Sardari}.
\newblock {Complexity of strong approximation on the sphere}.
\newblock {\em ArXiv e-prints}, March 2017.

\bibitem[{Sar}18]{least}
N.~T. {Sardari}.
\newblock {The least prime number represented by a binary quadratic form}.
\newblock {\em ArXiv e-prints}, March 2018.

\bibitem[Sel43]{Selberg}
Atle Selberg.
\newblock On the normal density of primes in small intervals, and the
  difference between consecutive primes.
\newblock {\em Arch. Math. Naturvid.}, 47(6):87--105, 1943.

\bibitem[Sou07]{Soundist}
K.~Soundararajan.
\newblock The distribution of prime numbers.
\newblock In {\em Equidistribution in number theory, an introduction}, volume
  237 of {\em NATO Sci. Ser. II Math. Phys. Chem.}, pages 59--83. Springer,
  Dordrecht, 2007.

\bibitem[TZ17]{Thorner}
Jesse Thorner and Asif Zaman.
\newblock An explicit bound for the least prime ideal in the {C}hebotarev
  density theorem.
\newblock {\em Algebra Number Theory}, 11(5):1135--1197, 2017.

\end{thebibliography}

\end{document}